\documentclass[a4paper,12pt]{article}

\usepackage{graphicx} 
\usepackage{amsmath} 
\usepackage{amssymb}  

\usepackage{comment}

\usepackage{bm} 

\usepackage{float}
\usepackage{mathtools}
\usepackage{algorithm}
\usepackage{algorithmic}

\usepackage{amsthm} 
\theoremstyle{definition}

\providecommand{\norm}[1]{\lVert#1\rVert}

\newcommand*{\defeq}{\mathrel{\vcenter{\baselineskip0.5ex \lineskiplimit0pt 
                     \hbox{\scriptsize.}\hbox{\scriptsize.}}}%
                     =} 
\newcommand{\such}{\,|\,}
\newcommand{\indicator}{1{\hskip -2.5 pt}\hbox{I}}

\title{\LARGE \bf
Measures and LMIs for optimal control of piecewise-affine systems
}

\begin{document}

\author{M. Rasheed Abdalmoaty$^1$,  Didier Henrion$^2$, Luis Rodrigues$^3$}

\footnotetext[1]{M. R. Abdalmoaty is with EADS Astrium GmbH; AOCS/GNC \& Flight Dynamics; 88039 Friedrichshafen, Germany. Most of this work was done when M. R. Abdalmoaty was with the Faculty of Electrical
Engineering, Czech Technical University in Prague, Czech Republic.}
\footnotetext[2]{D. Henrion is with CNRS, LAAS, 7 avenue du colonel Roche, F-31400 Toulouse, France;
Univ. de Toulouse, LAAS, F-31400 Toulouse, France. He is also with the Faculty of Electrical
Engineering, Czech Technical University in Prague,
Technick\'a 2, CZ-16626 Prague, Czech Republic. Corresponding author e-mail: {\tt henrion@laas.fr}.}
\footnotetext[3]{L. Rodrigues is with the Department of Electrical and Computer Engineering;
Concordia University; Montr\'eal, QC, H3G 2W1, Canada.}

\date{ \today}

\maketitle

\begin{abstract}
This paper considers the class of deterministic continuous-time optimal control problems (OCPs) with piecewise-affine (PWA) vector field, polynomial Lagrangian and semialgebraic input and state constraints. The OCP is first relaxed as an infinite-dimensional linear program (LP) over a space of occupation measures. This LP, a particular instance of the generalized moment problem, is then approached by an asymptotically converging hierarchy of linear matrix inequality (LMI) relaxations. The relaxed dual of the original LP returns a  polynomial approximation of the value function that solves the Hamilton-Jacobi-Bellman (HJB) equation of the OCP. Based on this polynomial approximation, a suboptimal policy is developed to construct a state feedback in a sample-and-hold manner. The results show that the suboptimal policy succeeds in providing a stabilizing suboptimal state feedback law that drives the system relatively close to the optimal trajectories and respects the given constraints. 
\end{abstract}

\section{Introduction}

Piecewise-affine (PWA) systems are a large modeling class for nonlinear systems.  However most of the nonlinear control theory does not apply to PWA systems because it requires certain smoothness assumptions. On the other hand, linear control theory cannot be simply employed due to the special properties of PWA systems inherited from nonlinear systems. PWA systems can naturally arise from linear systems in the presence of state saturation or from simple hybrid systems with state-based switching where the continuous dynamics in each regime are linear or affine \cite{Johansson}. Many engineering systems fall in this category, like power electronics converters for example. In addition, common electrical circuits components as diodes and transistors are naturally modeled as piecewise-linear elements.
PWA systems are also used to approximate large classes of nonlinear systems as in \cite{pedro}, \cite{luisphd}, and \cite{Casselman}. These approximations are then used to pose the controller design problem of the original nonlinear system as a robust control problem of an uncertain nonlinear system as suggested in \cite{Samadi,Samadi2}.

Problems of piecewise-affine systems are known to be challenging. The problems have a complex structure of regions stacked together in the state-space with each region containing an affine system. Any approach must identify the behavior in each region and then link them together to form a global picture of the dynamics. In \cite{blondel}, it has been shown that even for some simple PWA systems the problem of analysis or design can be either NP-hard or undecidable.

The motivation behind the research in this paper is the need of optimal control synthesis methods for continuous-time PWA systems with both input and state constraints. Additionally, there is often a need to find suitable tools for the design of stabilizing feedback controllers under state and input constraints. Over the last few years, there were several attempts addressing the synthesis problem for continuous-time PWA systems. The techniques are based on analysis methods and use convex optimization. These methods result in a state-based switched linear or affine controllers. For example, in \cite{boyd} a piecewise-linear state-feedback controller synthesis is done for piecewise-linear systems by solving convex optimization problem involving LMIs. The method is based on constructing  a globally quadratic Lyapunov function such that the closed-loop system is stable. Similarly, in \cite{rantzer} a quadratic performance index is suggested to obtain lower and upper bounds for the optimal cost using any stabilizing controller. However, the optimal controller is not computed. The method assumes a piecewise-affine controller structure which can be shown not to be always optimal (see section \ref{example}). In \cite{luis-how}, the work done in  \cite{boyd} is extended  to obtain dynamic output feedback stabilizing controllers for piecewise-affine systems. It formulates the search for a piecewise-quadratic Lyapunov function and a piecewise-affine controller as a nonconvex Bilinear Matrix Inequality (BMI), which is solved only locally by convex optimization methods. More recently in \cite{kamri}, a nonconvex BMI formulation is used to compute a state feedback control law.

For constrained PWA systems in discrete-time where both the partitions and the constraints are polyhedral regions, \cite{baotic} combines multi-parametric programming, dynamic programming and polyhedral manipulation to solve optimal control problems for linear or quadratic performance indices\footnote{We are grateful to Michal Kvasnica for pointing out this reference.}. The resulting solution when applied in receding horizon fashion guarantees stability of the closed-loop system.

To the best of the authors' knowledge there are no available guaranteed methods for synthesis of optimal controllers in continuous-time for PWA systems that consider general semi-algebraic state-space partitions, or that do not restrict the controller to be piecewise-affine, or do not require the performance index to be quadratic or piecewise-quadratic. 

The technique presented in this paper provides a systematic approach, inspired by \cite{lasserre:1643,carlo}, to synthesize a suboptimal state feedback control law for continuous-time PWA systems with multiple equilibria based on a polynomial approximation of the value function. The OCP is first formulated as an infinite-dimensional linear program (LP) over a space of occupation measures. The PWA structure of the dynamics and the state-space partition are then used to decompose the occupation measure of the trajectory into a combination of local occupation measures, one measure for each partition cell. Then, the LP formulation can be written in terms of only the moments of the occupation measures (countably many variables). This allows for a numerical solution via a hierarchy of convex LMI relaxations with vanishing conservatism which can be solved using off-the-shelf SDP solvers. The relaxations give an increasing sequence of lower bounds on the optimal value. An important feature of the approach is that state constraints  as well as any input constraints are very easy to handle. They are simply reflected into constraints on the supports of the occupation measures. It turns out that the dual formulation of the original infinite-dimensional LP problem on occupation measures can be written in terms of Sum-of-Squares (SOS) polynomials, that when solved yields a polynomial subsolution of the Hamilton-Jacobi-Bellman (HJB) equation of the OCP. This gives a good polynomial approximating value function along optimal trajectories that can be used to synthesize a suboptimal, yet admissible, control law. The idea behind the developed suboptimal policy exploits the structure of the HJB equation to generate the optimal control trajectory. The right-hand side of the HJB equation is iteratively minimized to construct a state feedback in a sample-and-hold manner with stabilization and suboptimality guarantees.


\section{The piecewise-affine optimal control problem}
\label{sec:pwaocp}

In this section, we first introduce the piecewise-affine continuous-time model, and then formulate the optimal control problem using some important assumptions.

\subsection{Piecewise-affine systems}
We consider exclusively continuous-time PWA systems. The term PWA is to be understood as PWA in the system state $x$. The state-space is assumed to be partitioned into a number of cells $X_i$ such that the dynamics in each cell  takes the form
\begin{equation}
\dot{x} = A_i x + a_i + B_i u \quad\quad \text{for} \quad  x \in X_i, \quad i \in I
\label{eq:pwa}
\end{equation}
where $I$ is the set of cell indices, and the union of all cells is $X=\cup_{ i \in I} X_i \subset \mathbb{R}^n$. The global dynamics of the system depends on both the cells and the corresponding local dynamics. The matrices $A_i$, $a_i$, and $B_i$ are time independent. In general, the geometry of the partition ${X_i}$ can be arbitrary. However, to arrive at useful results we assume the cells to be compact basic semi-algebraic sets (intersection of polynomial sublevel sets) with disjoint interiors. They are allowed to share boundaries as long as these boundaries have Lebesgue measure zero in $X$. There are many notions of solutions for PWA systems with different regularity assumptions on the vector field. The concern here is to ensure the uniqueness of the trajectories. Systems with discontinuous right-hand sides can have attracting sliding modes, non-unique trajectories, or trajectories may not even exist in the classical sense \cite{Cortes}, \cite{Johansson}. We assume that the PWA system is well-posed in the sense that it generates a unique trajectory for any given initial state. This is guaranteed if we assume that the global vector field is Lipschitz. It is usually the case if the model is the result of approximating a nonlinear function. 

\subsection{ Problem formulation}
 Optimal control problems (OCPs) of PWA systems are usually Lagrange problems where the state-space is partitioned into a finite number of cells. In addition to the dynamics, the cost functional can be also defined locally in each cell. The objective is to find optimal trajectories starting from an initial set and terminating at a target set that minimize the running cost and respect some input and state constraints. Consider the following general free-time PWA OCP with both terminal and running costs
\begin{equation}
\begin{aligned}
v^\ast(x_0) = &\underset{T,\: u}{\text{inf}}
& & L_T(x_T) +  \sum_{i=1}^r \int_{0}^T L_i(x(t),u(t))\indicator_{X_i}(x(t))dt  \\
& \text{s.t.}
& &  \dot{x} = A_i x(t) + a_i + B_i u(t), \quad x(t) \in X_i\\
& & & i=1,\ldots,r \quad t \in [0,T] \\
& & &x(0) = x_0 \in X_0 \subset \mathbb{R}^n,  & \\
& & & x(T) =x_T \in X_T \subset \mathbb{R}^n, \\
& & & (x(t), u(t)) \in X \times U \subset \mathbb{R}^n \times \mathbb{R}^m .
\label{pwaocp}
\end{aligned}
\end{equation}
where the scalar mapping  $\indicator_{X_i}: {X_i} \to \{0,1\}$ is the indicator function defined as follows:
 \begin{equation*}
\indicator_{X_i }(x(t)) \defeq \left\{ 
  \begin{array}{l l}
     1, & \quad \text{if} \quad x(t) \in X_i\\
     0, & \quad \text{if} \quad x(t) \notin X_i.\\
   \end{array} \right.
 \end{equation*}

The infimum in (\ref{pwaocp}) is sought over all admissible control functions $u(\cdot)$ with free final time $T$.
The dynamic programming approach of optimal control reduces the above problem to the problem of solving the following system of Hamilton-Jacobi-Bellman (HJB) equations for the value function $v^\ast$:

 \begin{equation}
 \begin{aligned}
  \inf_{u \in U} \{\nabla v^\ast(x) \cdot f_i(x,u) + L_i(x,u)\} = 0  \\
  \forall  (x,u) \in X_i \times U, \quad \forall i=1,\ldots,r\\
\end{aligned}
\label{hjb}
\end{equation}
with the terminal condition
\begin{equation*}
v^\ast(x(T)) =  L_T(x^\ast(T)).
\end{equation*}

In full generality, solving the HJB equations is very hard and the value function is not necessarily differentiable. Therefore, solutions must be interpreted in a generalized sense.

To proceed, we adopt  the following assumptions:
\begin{itemize}
        \item The terminal time $T$ is finite and the control functions $u(\cdot)$ are measurable.
        \item The PWA system is well-posed in the sense that it generates a unique trajectory from every initial state, and the vector field is locally Lipschitz.
        \item The Lagrangians  and the terminal cost are polynomial maps, namely $L_i \in \mathbb{R}[x,u]\; \forall i$ and $L_T \in     \mathbb{R}[x]$. 
                        \item The cells $X_i$, the control set $U$, the sets $X_0$ and $X_T$ are compact basic semi-algebraic sets defined as follows:
                                \begin{equation}\label{regions}
                         \begin{aligned}
X_i \times U = \{ (x,u)& \in \mathbb{R}^n  \times \mathbb{R}^m  \,|\, p_{i,k}(x,u) \geq 0,\\
&   \forall k = 1, \ldots, m_i\}, \; i=1,\ldots,r
                \end{aligned}
                \end{equation}
                
                          \begin{equation}
                 \begin{aligned}
                 X_0 &= \{ x \in \mathbb{R}^n \,|\, p_{0,k}(x) \geq 0,  \forall k = 1, \ldots, m_0 \},\\
                X_T &= \{ x \in \mathbb{R}^n \,|\, p_{T,k}(x) \geq 0,  \forall k = 1, \ldots, m_T \}.\\
                \end{aligned}
                \label{eq:sets}
                \end{equation}
         \item The cells $X_i$ have disjoint interior and they are allowed to share boundaries as long as these boundaries have Lebesgue measure zero in $X$.
   \end{itemize}

\section{The moment approach}
In this section, we formulate the nonlinear and nonconvex PWA OCP (\ref{pwaocp}) into a convex infinite-dimensional optimization problem over the state-action occupation measure $\mu$. The problem is then approached by an asymptotically converging hierarchy of LMI relaxations to arrive at a polynomial approximation of the value function.

\subsection{Occupation measures}
\label{subsec:}
Occupation measures are used to deal with dynamic objects where time is involved. We focus on the application of occupation measures to dynamic control systems where the dynamic objects are ordinary differential equations. Starting by the PWA vector field, we generate a sequence of moments by writing the ODEs in terms of occupation measures. We then manipulate the measures through their moments to optimize over system trajectories using a linear matrix inequality (LMI) representation.

To illustrate the measures formulation, first consider the uncontrolled autonomous dynamic system defined by the Lipschitz vector field $f:\mathbb{R}^n \to \mathbb{R}^n$ and the nonlinear differential equation
\begin{equation}
\dot{x} = f(x), \quad \quad x(0) = x_0.
\label{eq:ode}
\end{equation}

We think of the initial state $x_0$ as a random variable in $\mathbb{R}^n$ modeled by a probability measure $\mu_0$ supported on the compact set $X_0$, i.e. a nonnegative measure $\mu_0$ such that $\mu_0(X_0) = 1$.  Then at each time instant $t$,  the state $x(t)$ can also be seen as a random variable ruled by a nonnegative probability measure $\mu_t$. Solving (\ref{eq:ode}) for the state trajectory $x(t)$ yields a family of trajectories starting in $X_0$ and ending at a final set $X_T$. 
The measure $\mu_t$ of a set can then be thought of the ratio of the volume of trajectory points that lie inside that set at time $t$ to the total volume of points at time $t$.
In particular, if the number of trajectory points one is considering is finite, then the measure $\mu_t$ of a set is nothing else than the ratio of trajectory points that lie in the set at time $t$ to the total number of trajectory points.
The family of measures $\mu_t$ can thus be thought of as a density of trajectory points and satisfies the following  linear first-order continuity PDE in the nonnegative probability measures space

\begin{equation}\label{transport}
\frac{\partial \mu_t}{\partial t} + \nabla \cdot (f\mu_t) = 0.
\end{equation}

The above equation is known as Liouville's equation or advection PDE, in which $\nabla \cdot (f\mu_t)$ denotes the divergence of measure $f \mu_t$ in the sense of distributions. It describes the linear transport of measures from the initial set to the terminal set.

The occupation measure of the solution over some subset $\mathcal{T} \times \mathcal{X}$ in the $\sigma$-algebra of $[0,T] \times X$  is simply defined  as the time integration of $\mu_t$ as follows
\begin{equation}\label{mudef}
\mu (\mathcal{T}\times\mathcal{X}) \defeq \int_{\mathcal{T}} \mu_t(\mathcal{X}) dt=\int_{\mathcal{T}}\indicator_{\mathcal X}(x(t))dt,
\end{equation}
where the last equality is valid when $x(t)$ is a single solution to equation (\ref{eq:ode}).
It is important to note that when the initial condition $x_0$ is deterministic, the occupation measure $\mu$ is the time spent by the solution $x(t)$ in the subset $\mathcal{X}$ when $t \in \mathcal{T}$. 
We see that the occupation measure can indicate when the solution is within a given subset. 


\subsection{The primal formulation}
\label{subsec:primal}

In the sigma-algebra of Borel sets (smallest sigma-algebra that contains the open sets) let $\mathcal{M}(X)$ denote the space of signed Borel measures supported on a compact subset $X$ of the Euclidean space. Furthermore let $\mathcal{C}(X)$ be the space of bounded continuous functions on $X$, equipped with the supremum norm.  Then, the space $\mathcal{M}(X)$ is the dual space of $\mathcal{C}(X)$ with the following duality bracket
\begin{equation}
\langle v,\mu \rangle = \int_X v d\mu, \quad \quad \forall (v, \mu) \in \mathcal{C}(X) \times \mathcal{M}(X). 
\end{equation}

Now consider the piecewise-affine optimal control problem (\ref{pwaocp}), and assume that the control trajectory $u(t)$ is admissible such that all the constraints of the OCP are respected. Accordingly, we define the state-action  local occupation measure (including sets on the space of control inputs) associated with the cell $X_i$  to be
\begin{equation}
\mu_i ([0,T] \times X_i \times U) = \int_{0}^{T}  \indicator_{X_i \times U} (x(t),u(t)) dt.
\end{equation}

The local occupation measure $\mu_i (X_i \times U)$ encodes the trajectories in the sense that it measures the total time spent by the trajectories $(x(t), u(t))$ in the admissible set $X_i \times U$. Furthermore, define the global state-action occupation measure for the trajectory $(x(t),u(t))$ as a linear combination of local occupation measures such that
\begin{equation}
\begin{aligned}
\mu([0,T] \times X \times U) &= \int_{0}^{T}  \indicator_{X \times U} (x(t),u(t)) dt\\
& = \sum_{i=1}^r \mu_i(X_i \times U).
\end{aligned}
\end{equation}
The indicator function $\indicator_{X \times U} (x(t),u(t))$ equals $1$ through the interval $[0,T]$.
The occupation measure of the whole state space is given by the terminal time $T$. 
In general we should make sure that $T$ is finite, otherwise the occupation measure of the whole state space 
may escape to infinity.
The initial and terminal occupation measures are probability measures supported on $X_0$ and $X_T$, respectively. 
 
The objective function of problem (\ref{pwaocp}) can now be rewritten in terms of these measures to get the following linear cost
\begin{equation}
J(x_0,u(t)) =  \sum_{i=1}^r \int_{X_i \times U} L_i d\mu_i+\int_{X_T} L_T d\mu_T.
\end{equation}
With duality brackets, it reads
\begin{equation}
J(x_0,u(t)) =  \sum_{i=1}^r \langle L_i,\mu_i \rangle + \langle L_T, \mu_T  \rangle. 
\end{equation}
If the Lagrangian is the same for all the cells, say $L$, the performance measure can be written in terms of the global state-action occupation measure as
\begin{equation}
J(x_0,u(t)) =  \int_{X \times U} L d\mu + \int_{X_T} L_T d\mu_T.
\end{equation}
Our next step is to determine the measure transport equation that encodes the PWA dynamics in the measure space. To do that, we define a compactly supported global test function $v \in \mathcal{C}^1(X)$.
Then for $i=1,\ldots,r$ we define a linear map $F_i: \mathcal{C}^1(X_i) \to \mathcal{C}(X_i \times U)$
\begin{equation}
F_i(v) \defeq -\dot{v} = - \nabla v \cdot (A_ix+a_i+B_iu).
\end{equation}
Integration gives the following
\begin{equation}
\begin{aligned}
\int_0^T dv                                                                                                    
                    &=  -\sum_{i=1}^r \int_{X_i \times U} F_i(v) d\mu_i\\
                    &= \int_{X_T} v d\mu_T  -  \int_{X_0} v d\mu_0.  
\end{aligned}
\end{equation}
Equivalently, we can write
\begin{equation}
\begin{aligned}
\sum_{i=1}^r \langle F_i(v), \mu_i \rangle + \langle v, \mu_T \rangle &=   \langle v, \mu_0 \rangle \forall v \in \mathcal{C}^1(X).
\end{aligned}
\end{equation}

Now define the following ($r+1$)-tuple, in which we gather all the local occupation measures as the first $r$ elements, and the terminal occupation measure as the last element
 \begin{equation*}
 \nu \defeq (\mu_1,\ldots,\mu_r,\mu_T).
 \end{equation*}
The piecewise-affine optimal control problem (\ref{pwaocp}) is equivalent to the following infinite-dimensional linear optimization problem over occupation measures:
\begin{equation}
\begin{aligned}
p^\ast =  \quad &\underset{\nu} {\text{inf}}
& &   \sum_{i=1}^r \langle L_i,\mu_i \rangle + \langle L_T, \mu_T  \rangle \\
& \text{s.t.} 
& &   \sum_{i=1}^r \langle F_i(v), \mu_i \rangle + \langle v, \mu_T \rangle &=   \langle v, \mu_0 \rangle \\
& & & \quad \quad \quad \forall v \in \mathcal{C}^1(X).
\label{eq:LPv}
\end{aligned}
\end{equation}
Furthermore, defining the linear mapping $\mathcal{L}: \mathcal{C}^1(X) \to \prod_{i=1}^r\mathcal{C}(X_i) \times \mathcal{C}(X_T)$, as $\mathcal{L}(v)=(F_1(v), \ldots, F_r(v), v)$ we can rewrite the constraint in (\ref{eq:LPv}) as
 \begin{equation}
\begin{aligned}
\langle (F_1(v), \ldots, F_r(v), v),\nu \rangle  &=   \langle \mathcal{L}(v),  \nu \rangle\\
&= \langle v,  \mathcal{L}^\ast(\nu) \rangle = \langle v, \mu_0  \rangle, \\
& \quad  \forall v \in \mathcal{C}^1(X).
\end{aligned}
\end{equation}
This defines the adjoint map $\mathcal{L}^\ast: \prod_{i=1}^r\mathcal{M}(X_i) \times \mathcal{M}(X_T) \to \mathcal{M}^\ast(X)$.
The measure transport equation (\ref{transport}) is then given by \cite{evans}
\begin{equation}
\begin{aligned}
\mathcal{L}^\ast(\nu) = \mu_0 = \sum_{i=1}^r \nabla \cdot (f_i \mu_i)  +  \mu_T 
\end{aligned}
\end{equation}
with the symbol $f_i$ denoting the dynamics in the cell with index $i=1,\ldots,r$.

Finally, define the tuple $c \defeq  (L_1, \ldots, L_r, L_T)$, and associate the piecewise-affine optimal control problem (\ref{pwaocp})  to the following infinite-dimensional linear program 
\begin{equation}
\begin{aligned}
p^\ast =  \quad &\underset{\nu}{\text{inf}}
& &   \langle c,\nu \rangle\\
& \text{s.t.} 
& &   \mathcal{L}^\ast(\nu) = \mu_0 \\
&&& \nu \succeq 0.
\label{primalLP}
\end{aligned}
\end{equation}
This is the primal formulation of the OCP in terms of  occupation measures of the trajectory $(x(t),u(t))$. The nonlinear nonconvex PWA OCP (\ref{pwaocp}) is therefore reformulated as an infinite-dimensional LP in the measure space.
 
\subsection{Moments and LMI relaxations}
\label{subsec:lmi}

The moments of the occupation measure $\mu$ are defined  by integration of monomials with respect to $\mu$ . The $\alpha$-th moment of $\mu$ over the support $X$ is given by
\begin{equation}
y_\alpha = \int_X x^\alpha d\mu, \quad \forall \alpha \in \mathbb{N}^n.
\label{moments}
 \end{equation} 
where $x^\alpha=\prod_{i=1}^n x_i^{\alpha_i}$.
Noting that $d\mu(x)=\mu(dx)$ and using the definition (\ref{mudef}) of $\mu$, moments can be rewritten as
\begin{equation*}
y_\alpha = \int_0^T [x(t)]^\alpha dt, \quad \forall \alpha \in \mathbb{N}^n
 \end{equation*} 
where $x(t)$ denotes the solution of the ODE starting at $x_0$. This follows from the definition of $\mu$. Therefore, if we can find the moments and handle the representation conditions (\ref{moments}), solving the moments gives the solution of the ODE because the infinite (but countable) number of moments uniquely characterize a measure (on a
compact set).

Note that LPs (\ref{eq:LPv}) and (\ref{primalLP}) are equivalent. To proceed numerically we restrict the continuously differentiable functions to be polynomial functions of the state. In other words, we consider $v \in \mathbb{R}[x] \subset \mathcal{C}^1(X)$. By this restriction, we obtain an instance of the generalized moment problem (GMP), i.e. an infinite-dimensional linear program over moments sequences corresponding to the occupation measures. It turns out \cite[Ch3]{{mpp}} that if the supports of the measures are compact basic semi-algebraic sets, the GMP can be approached using an asymptotically converging hierarchy of LMI relaxations. 

To write the semidefinite relaxation of the primal infinite-dimensional LP, let $y_i = (y_{i_\alpha})$, \mbox{$\alpha \in  \mathbb{N}^n \times \mathbb{N}^m$} be the moments sequence corresponding to the local occupation measure $\mu_i$, $i=1,\ldots,r$. Moreover, let $y_0 = (y_{0_{\beta}})$ and $y_T = (y_{T_{\beta}})$ with $\beta \in \mathbb{N}^n$ be the moment sequences corresponding to $\mu_0$ and $\mu_T$ respectively. Given any infinite sequence $y=(y_\alpha)$ of real numbers with $\alpha \in \mathbb{N}^n$, 
define the linear functional $ \ell: \mathbb{R}[x] \to \mathbb{R}$ that maps polynomials to real numbers as follows
\begin{equation*}
p(x) = \sum_{\alpha \in \mathbb{N}^n} p_\alpha x^\alpha \quad \mapsto \quad \ell_y(p) = \sum_ {\alpha \in \mathbb{N}^n} p_\alpha y_\alpha .
\end{equation*}

In each LMI relaxation we truncate the infinite moment sequence to a finite number of moments. The LMI relaxation of order $d$, including moments up to $2d$, of the GMP instance (\ref{primalLP}) can be formulated by taking test functions $v = x^\alpha$ with $\alpha \in \mathbb{N}^n$, such that $\text{deg } v = 2d$, as follows
\begin{equation}
\begin{aligned}
p^\ast_d =   &\underset{y_1,\ldots,y_r,y_T}{\text{inf}}
& &   \sum_{i=1}^r \ell_{y_i}(L_i) + \ell_{y_T}(L_T)\\
& \text{s.t.}
& &  \sum_{i=1}^r \ell_{y_i}(F_i(v)) + \ell_{y_T}(v) = \ell_{y_0}(v), \\
& & & M_d(y_i) \succeq 0,  \forall i\\
& & & M_d(p_{i,k} \: y_i) \succeq 0, \forall i,  \forall k = 1, \dots, m_g,\\
& & & M_d(y_T) \succeq 0,\\
& & &  M_d(p_{T,k}\: y_T) \succeq 0, \forall  k = 1, \dots, m_T.
\end{aligned}
\label{primalSDP}
\end{equation}

The minimum relaxation order has to allow the enumeration of all the moments appearing in the objective function and the linear equality constraint. The matrices $M_d(y_i)$ and $M_d(y_T)$ are called moment matrices of the local occupation measure and terminal probability measures, respectively. Each moment matrix is defined to be a square matrix of dimension ${{d+n}\choose{n}}$ filled with the first $2d$ moments corresponding to the representing measure. They are linear in the moments. Similarly, the matrices $M_d(p_{i,k} \: y_i)$ and $M_d(p_{T,k}\: y_T)$ are linear in the moments and are called localizing matrices. The linear equality constraint represents the peicewise-affine dynamics, and the LMIs ensure that (\ref{moments}) holds and that the measures are supported on the given sets defined in (\ref{regions})--(\ref{eq:sets}) (see \cite{mpp} for more details). 

\subsection{The dual formulation}
\label{subsec:dual}
The duality between finite measures and compactly supported  bounded continuous functions is captured by convex analysis. The dual of the LP (\ref{primalLP}) is thus formulated over the space of positive bounded continuously differentiable functions as follows
\begin{equation}
\begin{aligned}
d^\ast =  \quad &\underset{v \in \mathcal{C}^1(X)}{\text{sup}}
& &   \langle v,\mu_0 \rangle\\
& \text{s.t.}
& & z = c - \mathcal{L}(v)\\
& & & z \geq 0
\label{dualLP}
\end{aligned}
\end{equation}
where $z$ is a vector of continuous functions.
The dual LP can be written in more explicit form to reveal the structure of the linear constraints. We can equivalently write 
\begin{equation}
\begin{aligned}
d^\ast =  \quad &\underset{v \in \mathcal{C}^1(X)}{\text{sup}}
& &   \langle v,\mu_0 \rangle\\
& \text{s.t.}
& &   L_i - F_i(v) \geq 0, \quad \forall i=1,\ldots,r,\\
& & & L_T - v(x_T) \geq 0
\end{aligned}
\end{equation}
and more explicitly
\begin{equation}\label{SDP}
\begin{aligned}
d^\ast =  \quad &\underset{v \in \mathcal{C}^1(X)}{\text{sup}}
& &   \int_{X_0} v d\mu_0\\
& \text{s.t.}
& & \nabla v(x)\cdot f_i + L_i(x,u) \geq 0,\\
& & & \quad \forall (x,u) \in X_i \times U,  \forall i=1,\ldots,r\\
& & & L_T - v(x_T) \geq 0,   \forall x \in X_T.
\end{aligned}
\end{equation}
We note that any feasible solution of SDP (\ref{SDP}) is actually a global smooth subsolution of the HJB equations (\ref{hjb}). 
By conic complementarity, along the optimal trajectory $(x^\ast,u^\ast)$, it holds
\begin{equation}
\langle  z^\ast, \nu^\ast \rangle = 0.
\end{equation}
Therefore, for the optimal dual function $v^\ast$, the following holds:
\begin{equation}
\begin{aligned}
\nabla v^\ast(x^\ast)\cdot f_i + L_i(x^\ast,u^\ast) = 0&, \\
 \forall  (x^\ast,u^\ast) \in X_i \times U&,\\
  \forall i=1,\ldots,r&
 \end{aligned}
\label{eq:pwaH}
\end{equation} 
and in addition,
\begin{equation}
v^\ast(x^\ast(T)) =  L_T(x^\ast(T)).\\
 \label{eq:pwaT}
\end{equation}

This is an important result. It shows the following:
\begin{enumerate}
  \item With a careful look, we can easily identify what we have in  (\ref{eq:pwaH}) to be the HJB PDE of the PWA optimal control problem satisfied along optimal trajectories, with the terminal conditions given by (\ref{eq:pwaT}).
 
\item The optimal dual function $v^\ast(x)$  is equivalent to the value function of the optimal control problem, hence the notation. The maximizer function $v^\ast(x)$ of the dual infinite-dimensional LP in equation (\ref{dualLP}) solves, globally, the HJB equation of the PWA optimal control problem along optimal trajectories.

\end{enumerate}

The dual convex relaxation, dual of LMI (\ref{primalSDP}), is formulated over positive polynomials. Putinar's Positivstellensatz \cite{putinar} is used to enforce positiveness. Therefore, the unknown dual variables are the coefficients of the polynomial $v$ and several SOS polynomials that deal with the polynomial positivity conditions of the constraints.
The dual program can then be written as follows:
\begin{equation}
\begin{aligned}
d_{d}^\ast =  \quad & \underset{v_d,s}{\text{sup}}
& &  \;  \int_{X_0} v_d d\mu_0\\
& \text{s.t.}
& &   L_i - F_i(v_d) = s_{i,0}+ \sum_{k=1}^{m_i} p_{i,k}\: s_{i,k}  \\
& & &\forall (x,u) \in X_i \times U, \quad  \forall i=1,\ldots,r,\\
& & & L_T - v_d(x) = s_{T,0}+ \sum_{k=1}^{m_T} p_{T,k}\: s_{T,k}   \\
& & & \forall x \in X_T.
\end{aligned}
\end{equation}
in which the degree of $v_d$ is $d$. The polynomials $s_{i,0},\; s_{i,k},\; s_{T,0} \text{ and } s_{T,k}$ are positive. They are Putinar's SOS representations of the constraint polynomials \cite{putinar}.   The polynomials  $p_{i,k}$ define the set $X_i \times U$ and the polynomials  $p_{T,k}$ define the set $X_T$, see equations (\ref{regions})--(\ref{eq:sets}). 

Using weak-star compactness arguments similar to those in the proof of \cite[Theorem 2]{roa},
we can show that $p^\ast = d^\ast$. Moreover, the equality $p^\ast=v^\ast(x_0)$ holds if
we allow relaxed controls (enlarging the set of admissible control functions
to probability measures) and chattering phenomena in OCP (\ref{pwaocp}), see e.g.
the discussion in \cite[Section 3.2]{roa} and references therein.
Note that in \cite{lasserre:1643}, the identity  $p^\ast = d^\ast = v^\ast(x_0)$ was shown
under stronger convexity assumptions.

\section{Suboptimal control synthesis}
\label{subsec:policy}

Assume that the analytical value function $v^\ast(x)$ of a general optimal control problem is available by solving the HJB PDE.
The optimal feedback control function $k^\ast(x(t))$ can then be selected such that it generates an admissible optimal control trajectory $u^\ast(t)$ that satisfies the optimal necessary and sufficient conditions (\ref{eq:pwaH}). The resulting optimal feedback function $k^\ast(\cdot)$ generates the admissible optimal control trajectory starting from any initial value $x_0$. Therefore we get the solution of the OCP as a feedback strategy, namely, if the system is at state $x$, the control is adjusted to $k^\ast(x)$. This approach has the difficulty that even if the value function was smooth there would be in general no continuous optimal state feedback law that satisfies the optimality conditions for every state.

It is well-known from Brockett's existence theorem that if a dynamic control system with some vector field $f(x,u)$ admits a continuous stabilizing feedback law (taking the origin as equilibrium point), then for every $\delta > 0$, the set $f(B(0,\delta), U)$ is a neighborhood of the origin \cite{brockett}, where $B(0,\delta)$ is a ball of radius $\delta$ with a center at zero. One famous example where this condition fails is the nonholonomic integrator \cite{Cortes}. More recently it was shown that asymptotic controllability of a system is necessary and sufficient for state feedback stabilization without insisting on the continuity of the feedback. The synthesis of such discontinuous feedbacks was described in \cite{clarke}, together with a definition of a solution concept for an ODE with discontinuous dynamics, namely the sample-and-hold implementation. It turns out that this concept is very convenient for PWA systems. The first paper to deal with sampled-date PWA systems in the form used here can be found in reference \cite{Rodrigues}.

 The suboptimal trajectory of the control system is defined by a partition, call it $\pi$, of the time set $[0,T_{\pi}]$ as follows: define $\pi = \{t_j\}_{0 \leq j \leq p}$ with a given diameter $d(\pi) \defeq \sup_{0 \leq j \leq p-1}(t_j, t_{j+1})$ such that $0=t_0<t_1 \dots <  t_p= T_\pi$. The sequence $t_j$ for $ 1\leq j \leq p$ depends on the evolution of the trajectory and the diameter of the partition $d(\pi)$.
Starting at an initial point $t_j$, the suboptimal trajectory is the classical solution of the Lipschitz ODE 
\begin{equation*}
 \begin{aligned}
&\dot{x}(t) = A_ix(t) + a + B_i k^\ast(x_j), \quad x(t) \in X_i\\
  & x_j = x(t_j) , \quad   t_j \leq t \leq t_{j+1}
\end{aligned}
\end{equation*}
such that $t_{j+1}$ depends on the evolution of the trajectory.  

The generated suboptimal trajectory corresponds to a piecewise constant open-loop control (point-wise feedback) which has physical meaning. The optimal points of the trajectory are those points at the beginning of each subinterval at which the control is updated using the optimal control function. The suboptimal trajectory converges to the optimal trajectory when $d(\pi) \to 0$. This corresponds to increasing the sampling rate in the implementation such that the sampling time $T_s \to 0$. In this case the generated trajectory corresponds to a fixed optimal feedback. 
This can be shown as follows: take a partition $\pi$ with $d(\pi)=\delta$, we can use the mean value theorem to write
\begin{equation*}
\begin{aligned}
v^\ast(x_{j+1}) - v^\ast(x_j) &= [\nabla v^\ast(x(\tau)) \cdot f_i(x(\tau), u_j^\ast(x_j))] (t_{j+1} - t_j)\\
 \delta_j &= (t_{j+1} - t_j) <= \delta, \text{ and } \tau \in [0, \delta_j].
\end{aligned}
\end{equation*}
and we note that when $\delta_j \to 0$,
\begin{equation*}
\begin{aligned}
\dot{v^\ast}(x_j) &\to \nabla v^\ast(x(t_j) \cdot f_i(x(t_j), u_j^\ast(x_j)\\
 &= - L_i(x_j,u_j^\ast(x_j)).
\end{aligned}
\end{equation*}
and the algorithm converges to the optimal trajectories.
It is then clear that for every initial condition $x_0$, and some prescribed $g > 0, \epsilon > 0$; there exist some $\delta > 0$ and $T_\pi >0$ such that if $d(\pi) < \delta$ the generated suboptimal trajectories starting at $x_0$ satisfy
\begin{equation*}
\begin{aligned}
J(x_0, u(t)) - v^\ast(x_0) &< g  \quad \text{(suboptimality gap)}\\
\norm{(x(T_\pi) - x_T)} &\le \epsilon \quad \text{(tolerance). } 
\end{aligned}
\end{equation*}
 In the special case of having $\delta = 0$, we have no gap due to the algorithm and $T_\pi = T$.
 
Assuming that a polynomial approximation of the value function for a given PWA continuous system is available by solving the relaxed dual LMI, the proposed suboptimal feedback policy $[0,T) \times \mathbb{R}^n \to U$  is constructed using closed-loop sampling and Algorithm \ref{alg1}.
\begin{algorithm}                      
\caption{Algorithm for Suboptimal Synthesis}          
\label{alg1}                           
\begin{algorithmic}                    
\STATE {\bf given:} $d(\pi)$, $\epsilon$, $x_0\in X_i$, $x_T$, and poly. approx. of $v^\ast(\cdot)$
\STATE {\bf initialization:} $t = t_0 = 0$, $x(0) = x_0$, $j =0$. \COMMENT{first interval}
\vspace{0.2cm}
\WHILE{$\norm{x(t) - x_T} > \epsilon$}
\STATE solve the static polynomial optimization problem
\begin{equation*}
u_j^\ast(x_j) = \underset{u \in U}{\operatorname{argmin}} \;  \nabla v(x_j) \cdot f_i(x_j,u) + L_i(x_j, u)
\end{equation*}
\REPEAT
\vspace{0.2cm}
\STATE solve $\dot{x}(t) = A_i x(t) + a_i + B_i u_j^\ast(x_j)$, starting at $t_j$
\vspace{0.2cm}
\UNTIL{$x(t) \notin X_i$ or $t-t_j = d(\pi)$ }
\vspace{0.2cm}

\STATE set $j \leftarrow j+1$ \COMMENT{next interval}
\STATE set $t_j \leftarrow t$, and $x_j \leftarrow x(t)$. \COMMENT{initial values for next interval}
\STATE determine new region $X_k$ for $x(t)$ and set $i \leftarrow k$
\ENDWHILE
\vspace{0.2cm}

\STATE $p \leftarrow j$ \COMMENT{number of intervals in $\pi$}
\RETURN{$\pi = \{t_j\}_{0 \leq j \leq p}$ }, $x(t)$, and $u(t)$ with $0<t\le T_\pi$.
\end{algorithmic}
\end{algorithm}

\section{Numerical Example}\label{example}

We consider a first-order PWA dynamic system with two cells. The optimization is done over both the control action and time horizon (free time). The OCP is defined as follows:
\begin{equation}
\begin{aligned}
v^\ast(x_0) =  \quad &\underset{T,u}{\text{inf}}
& & \int_0^T (2(x-1)^2+u^2) dt \\
& \text{s.t.}
& &  \dot{x}  = \left\{ \begin{array}{ll}
                                   f_1 =   -x + 1 +u, &  x \in X_1\\
                                     f_2 = x + 1 +u,  &   x \in X_2\\
                                    \end{array} \right. \\
& & &x(0) = -1,\,\,\,\,x(T) =+1
\label{eq:ex1}
\end{aligned}
\end{equation}
where $v^\ast$ is the value function, and $x_0 = x(0)$ is the initial state. The control $u(t)$ is assumed to be a bounded measurable function defined on the interval $t \in [0,T]$, and taking its values in $\mathbb{R}$. The state space is partitioned into two unbounded regions
\begin{displaymath}
\begin{aligned}
X &= X_1 \cup X_2, \quad \text{where,}\\
X_1 &= \{ x \in X \such x \geq 0\} \text{,}\\
X_2 &= \{ x \in X \such x \leq 0\}.
\end{aligned}
\end{displaymath}

The control Hamiltonian
\begin{equation} \label{h_ex1}
 H(x,u,\nabla v^\ast) = \left\{  \begin{array}{l l l}
   & \nabla v^\ast f_1 + (2(x-1)^2+u^2) , x \in X_1\\
   & \nabla v^\ast f_2 + (2(x-1)^2+u^2),  x \in X_2\\
 \end{array}  \right.\\\end{equation}
is used to write the HJB equation. Since the terminal time $T$ is subject to optimization, the control Hamiltonian vanishes along the optimal trajectory.

This gives the HJB PDE
\begin{equation} \label{eq:hjb_ex1}
\underset{u \in \mathbb{R}}{\text{inf}} \, H(x(t),u(t),\nabla v^\ast(x(t))) = 0,  \; \forall x(t), \; t \in [0,T].
 \end{equation}


An analytical optimal solution can be obtained by solving the HJB equations corresponding to each cell. This results in a state feedback control
\begin{equation}\label{eq:analytical_optimal}
k^\ast(x) = \left\{  \begin{array}{l c}
                                      (1-\sqrt{3})(x-1) &  \text{if} \quad x \geq 0\\
                                     -x - 1 + \sqrt{2(x-1)^2 + (x+1)^2} &  \text{if} \quad x \leq 0\\
                                    \end{array}  \right.\\\end{equation}
that satisfies the HJB PDE (\ref{eq:hjb_ex1}) for all $t \in [0,T]$ and $x$. The optimal control trajectory
\begin{displaymath}
 u^\ast(t) = k^\ast(x^\ast(t)) \quad \quad \forall t \in [0,T]
\end{displaymath} 
is obtained when the optimal state feedback is applied starting at the given initial condition. 
 The partial derivative of the value function with respect to the state is a viscosity solution of the HJB PDE (\ref{eq:hjb_ex1}), see e.g. \cite{evans}.
                                    
The initial and final states are known. Therefore, the initial and terminal occupation measures are Dirac measures. The  global occupation measure $\mu  \in {\mathcal M}(X\times U)$ is used to encode the system trajectories. It is defined as a combination of two local occupation measures, one for each cell, as follows
\begin{displaymath}
\begin{aligned}
\mu &= \mu_1 + \mu_2
\end{aligned}
\end{displaymath}
such that the support of $\mu_1$ is $\{ (x,u) \such x \in X_1, u \in U\}$ and the support of $\mu_2$
is $\{ (x,u) \such x \in X_2, u \in U\}$.
The OCP can then be formulated as an infinite-dimensional LP in measure space
\begin{equation} \label{eq:lp1}
\begin{aligned}
&p^\ast = \underset{\mu_1,\mu_2}{\text{inf}} \int_{X_1,U}L d\mu_1 +  \int_{X_2,U} Ld\mu_2 \\
 &\text{s.t.}  \int \nabla v \cdot f_1 d\mu_1 +  \int \nabla v \cdot f_2 d\mu_2 =\\
 & \quad\quad\quad\quad\quad\quad \int_{X_T} v d\mu_T -\int_{X_0} v  d\mu_0, \quad \forall v
\end{aligned}
\end{equation}
where $\mu_0$ and $\mu_T$ are the initial and final measures, supported on $X_T = \{+1\}$, $X_0 = \{-1\}$, respectively, and $L=2(x-1)^2+u^2$ is the Lagrangian. The functions $v$ are functions of the state $x$ and belong to the space of continuously differentiable  functions. The Matlab toolbox GloptiPoly \cite{gloptipoly}
is used  to formulate this infinite-dimensional LP on measures as a GMP. Here, we are mainly interested in solving the dual formulation on functions.

In the present example, both the state space and the input space are not compact. Hence, the LMI relaxations do not include any localizing constraints; moreover, Putinar's conditions are not sufficient for the convergence of the LMI optimal value. The sufficiency is guaranteed only for measures on compact support. The main concern here is that the mass (zeroth moment) of the occupation measure is not bounded and might tend to infinity. This is numerically problematic. There are two ways to turn this around. The first is to make the time interval compact by adding a constraint on the mass of the global occupation measure. For example we enforce that the mass of $\mu$ is less than a given large positive number. This exploits the fact that the vector field is asymptotically stable. Another way, which is more general, is to constrain $X$ and $U$ to sufficiently large subsets. One possibility is to take $X$ as a large ball  centered around $0$ (or the initial state $x_0$). The same can be assumed for the input space $U$. This is done without any problems as long as the optimal trajectory remains in the constraint set $(X,U)$.
To avoid numerical problems, it is also recommended in all cases to scale down the problem, if possible, to have all the variables inside the unit box. Scaling avoids blow up of the moments sequences for high relaxation orders.

The results obtained below for the OCP (\ref{eq:ex1}) assume that the global occupation measure is supported on $X \times U$ without introducing any constraints on the spaces.

The solution consists of two main steps:
\begin{enumerate}
\item Finding a relatively good polynomial approximation of the value function by solving the dual LMI relaxation of the infinite-dimensional LP (\ref{eq:lp1}).
\item Employing the suboptimal strategy developed in section \ref{subsec:policy} to generate a suboptimal admissible feedback control law based on the smooth approximation of the value function obtained in step 1. 
\end{enumerate}

\begin{figure}[t]
\begin{center}
\includegraphics[width=.8\textwidth]{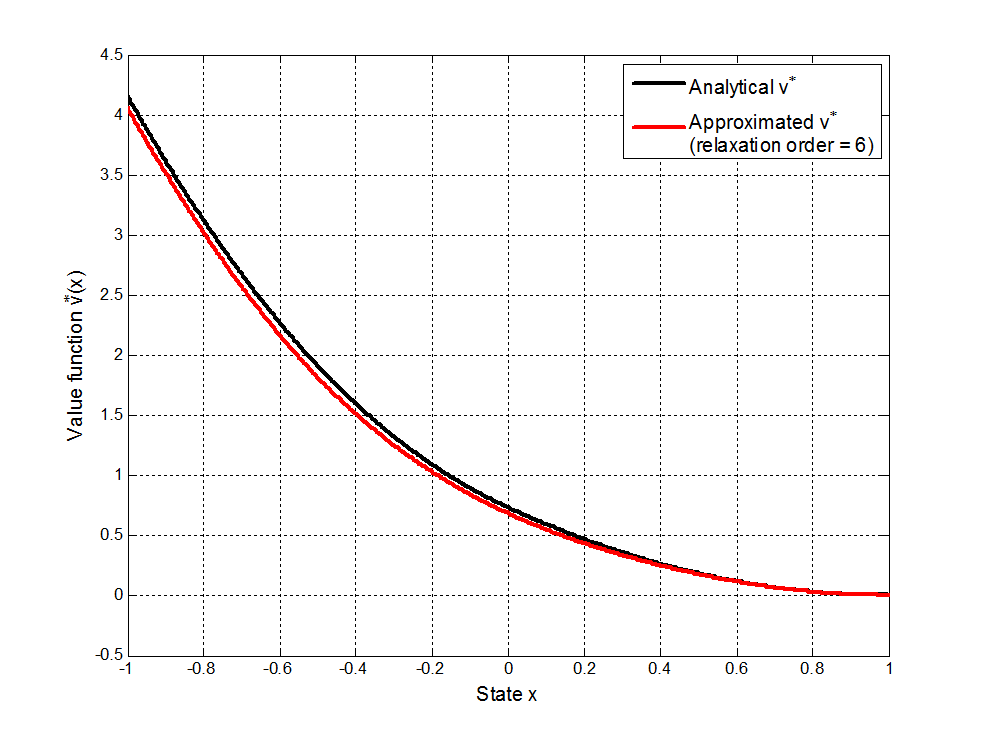}
\caption{PWA system. Approximating value function for $d = 6$.}
\label{fig:ex1_vf}
\end{center}
\end{figure}

The first step is achieved by solving the relaxed dual LMI on continuous functions. Figure \ref{fig:ex1_vf} shows the obtained approximation of the value function with LMI relaxation order $d=6$.

The obtained approximation is a polynomial of the state $x$ with degree equal to $2d$.

\begin{figure}[t]
\begin{center}
\includegraphics[width=0.8\textwidth]{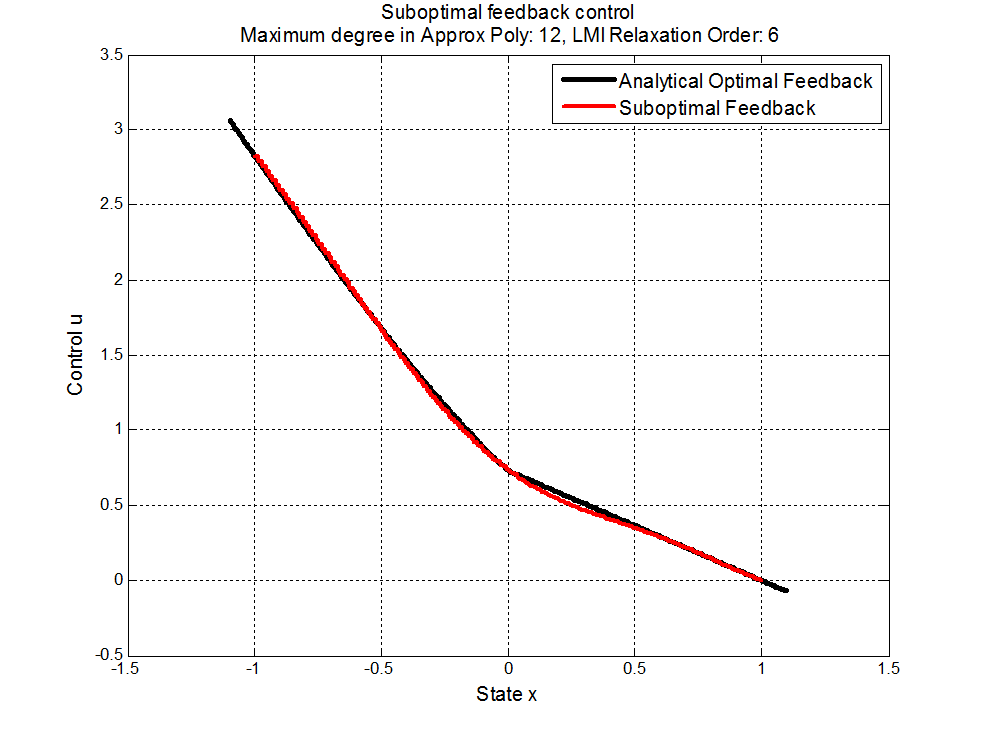}
\caption{PWA system. Suboptimal feedback for $d = 6$.}
\label{fig:ex1_fb}
\end{center}
\end{figure}

Based on this approximation, we employ algorithm \ref{alg1} from section \ref{subsec:policy}. The resulting suboptimal feedback control is shown in figure \ref{fig:ex1_fb} in comparison with the analytical optimal feedback, $k^\ast(x)$, calculated using  (\ref{eq:analytical_optimal}). It is clear that the algorithm gives a suboptimal feedback close to the optimal.

\section{Conclusions}

The focus of this paper is the synthesis of suboptimal state feedback controllers for continuous-time optimal control problems (OCP) with piecewise-affine (PWA) dynamics and piecewise polynomial cost functions. Both state constraints and input constraints are considered in a very convenient way and they do not pose additional complexity.

The problem is formulated as an abstract infinite-dimensional convex optimization problem over a space of occupation measures which is then solved via a converging hierarchy of LMI problems. By restricting the dual variables in the dual of the original infinite-dimensional program  to be monomials, we obtain a polynomial representation of the value function of the OCP in terms of upper envelope of subsolutions to a system of HJB equations corresponding to the OCP. By fixing the degree of the monomials, the same dual program can be relaxed and written, using Putinar's Positivstellensatz, as a polynomial sum-of-squares (SOS) program, which can be transformed and solved as an LMI problem.
As soon as the polynomial approximation of the value function is available, one can systematically generate a suboptimal, yet admissible, feedback control. The suboptimal control strategy is based on a closed-loop sampling implementation which is very convenient for PWA systems. The method generates an approximate control signal which is piecewise constant, and near optimal trajectories that respect the given constraints.

\section*{Acknowledgments}

Part of this work was
supported by project number 103/10/0628 of the Grant Agency of the Czech Republic,
and we are particularly grateful to Martin Hrom\v c\'{\i}k for his help.



\begin{thebibliography}{99}

%

\bibitem{baotic}
L. Baoti\'{c}.
Optimal control of piecewise-affine systems - a Multi-parametric Approach -
PhD thesis, ETH Zurich, 2005.

\bibitem{blondel}
V. D. Blondel , J. N. Tsitsiklis,
 Complexity of stability and controllability of elementary hybrid systems,
 Automatica, vol. 35, no. 3, pp. 479-489, 1999.
  



\bibitem{brockett}
R. W. Brockett.
Asymptotic stability and feedback stabilization.
In Differential Geometric Control Theory (R. W. Brockett, R. S. Millman and H. J. Sussmann, eds.),
pp 181-191, Birkhauser, Boston, MA, 1983.


\bibitem{Casselman}
S. Casselman, L. Rodrigues. 
A new methodology for piecewise affine models using Voronoi partitions.
Proceedings of the IEEE Conference on Decision and Control, December 2009.

\bibitem{clarke}
F. Clarke, Y. Ledyaev, E. D. Sontag, A. Subbotin.
Asymptotic controllability implies feedback stabilization.
IEEE Trans. Automat. Control, vol. 42, no. 10, pp. 1394-1407, 1997.

\bibitem{Cortes}
J. Cort{\'e}s.
Discontinuous dynamical systems: a tutorial on solutions, nonsmooth analysis, and stability.
IEEE Control Syst. Mag., vol.28, no. 3, pp. 36-73, 2008.
 
 
 
 
 
\bibitem{evans}
L. C. Evans.
Partial differential equations. 
2nd edition. Graduate Studies in Mathematics, vol. 19. AMS, Providence, RI, 2010.
 
\bibitem{boyd} 
A. Hassibi, S. Boyd.
Quadratic stabilization and control of piecewise-linear systems.
Proceedings of the American Control Conference, January 1998.
 

\bibitem{gloptipoly}
D. Henrion, J. B. Lasserre, J. L\"ofberg.
GloptiPoly 3: moments, optimization and semidefinite programming.
Optim. Methods and Software, 24(4-5):761-779, 2009.
 
\bibitem{Johansson}
M. Johansson.
Piecewise linear control systems. 
A computational approach. LNCIS 284. Springer Verlag, Berlin, 2003.

\bibitem{pedro}
P. Juli{\'a}n, A. Desages, O. Agamennoni.
High-level canonical piecewise linear representation using a simplicial partition.
IEEE Trans. Circuits Systems I Fund. Theory Appl., vol. 46, no. 4, pp. 463-480, 1999.

\bibitem{kamri}
D. Kamri, R. Bourdais, J. Buisson, C. Larbes.
 Practical stabilization for piecewise-affine systems: A BMI approach.
Nonlinear Analysis: Hybrid Systems, vol. 6, no. 3 Pages 859-870, 2012.

\bibitem{carlo}
D. Henrion, J. B. Lasserre, C. Savorgnan. Nonlinear optimal control synthesis
via occupation measures. Proc. IEEE Conf. Dec. Control, Cancun, Mexico, Dec. 2008. 

\bibitem{roa}
D. Henrion, M. Korda. Convex computation of the region of attraction of polynomial control systems.
LAAS-CNRS Research Report 12488, Sep. 2012.


\bibitem{mpp}
J. B. Lasserre. 
Moments, positive polynomials and their applications. 
Imperial College Press, London, UK, 2009.

\bibitem{lasserre:1643}
J. B. Lasserre, D. Henrion, C. Prieur, E. Tr\'elat.
Nonlinear optimal control via occupation measures and LMI relaxations.
SIAM J. Control Optim., vol. 47, no. 4, pp. 1643-1666, 2008.



\bibitem{putinar}
M. Putinar.
Positive polynomials on compact semi-algebraic sets. 
Indiana Univ. Math. J., vol. 42, no. 3, pp. 969-984, 1993.

\bibitem{rantzer}
A. Rantzer, M. Johansson.
Piecewise linear quadratic optimal control.
IEEE Trans. Automat. Control, vol. 45, no. 4, pp. 629-637, 2000.

\bibitem{luisphd}
L. Rodrigues.
Dynamic output feedback controller synthesis for piecewise-affine systems.
PhD thesis, Stanford Univ., Palo Alto, CA, 2002.

\bibitem{luis-how} 
L. Rodrigues, A. Hassibi, J. P. How. 
Output feedback controller synthesis for piecewise-affine systems with multiple equilibria.
Proceedings of the American Control Conference, June 2000.

\bibitem{Samadi}
B. Samadi, L. Rodrigues.
Extension of local linear controllers to global piecewise affine controllers for uncertain non-linear systems.
Int. J. Systems Sci., vol. 39, no. 9, pp. 867-879, 2008.

\bibitem{Samadi2}
B. Samadi and L. Rodrigues.
Controller Synthesis for Piecewise-Affine Slab Differential Inclusions: a Duality-Based Approach.
Automatica, 45 (7), pp.812-816, 2009.

\bibitem{Rodrigues}
L. Rodrigues.
Stability of sampled-data piecewise-affine systems under state feedback.
Automatica, (43) 7, pp.1249-1256, 2007.








\end{thebibliography}
\end{document}